\documentclass{amsart}
\usepackage{amsfonts}
\usepackage{amsmath}
\usepackage{amssymb}
\usepackage{graphicx}
\usepackage[all,dvips]{xy}

\newtheorem{thm}{Theorem}[section]
\newtheorem{cor}[thm]{Corollary}
\newtheorem{lem}[thm]{Lemma}
\newtheorem{prop}[thm]{Proposition}
\theoremstyle{definition}
\newtheorem{defn}[thm]{Definition}
\theoremstyle{remark}
\newtheorem{rem}[thm]{Remark}

\newtheorem{exs}[thm]{Examples}

\newenvironment{dem}{\noindent {\bf Proof.}}{\hfill $\Box$\\}

\numberwithin{equation}{section}

\newcommand{\Ann}{\mathrm{Ann}}
\newcommand{\add}{\mathrm{add}}
\newcommand{\md}{\mathrm{mod}}
\newcommand{\ind}{\mathrm{ind}}
\newcommand{\Hom}{\mathrm{Hom}}
\newcommand{\Ext}{\mathrm{Ext}}
\newcommand{\End}{\mathrm{End}}
\newcommand{\HH}{\mathrm{HH}}
\newcommand{\rad}{\mathrm{rad}}
\newcommand{\pd}{\mathrm{pd}}
\newcommand{\id}{\mathrm{id}}
\newcommand{\gd}{\mathrm{gl.dim}}
\newcommand{\LA}{\mathcal{L}_A}
\newcommand{\RA}{\mathcal{R}_A}
\newcommand{\LUR}{\mathcal{L}_A \cup \mathcal{R}_A}

\begin{document}

\sloppy

\title{Algebras determined by their supports}

\author  [I.~Assem]{Ibrahim Assem}
\address{%
Ibrahim Assem\\
D\'epartement de Math\'ematiques\\%
Universit\'e de Sherbrooke\\%
Sherbrooke, Qu\'ebec\\%
Canada, J1K 2R1.}
\email{ibrahim.assem@usherbrooke.ca}

\author[D.~Castonguay]{Diane Castonguay}
\address{%
Diane Castonguay\\ Instituto de Informática\\%
Universidade Federal de Goiás, Campus II - Samambaia, CEP: 74001-970\\%
Goi\^ania, Brazil}%
\email{diane@inf.ufg.br}

\author[M.~Lanzilotta]{Marcelo Lanzilotta}
\address{%
Marcelo Lanzilotta\\ Centro de Matem\'atica (CMAT)\\%
Igu\'a 4225\\ Universidad de la Rep\'ublica, CP 11400\\%
Montevideo, Uruguay.} %
\email{marclan@cmat.edu.uy}

\author[R.~Vargas]{Rosana R. S. Vargas}%
\address{%
Rosana Vargas\\ Escola de Artes, Ciências e Humanidades (EACH)\\%
Universidade de S\~ao Paulo, CEP: 03828-000\\%
S\~ao Paulo, Brazil}%
\email{rosanav@usp.br}

\thanks{
The authors wish to thank E. R. Alvares for fruitful discussions. The first author gratefully acknowledges partial support from the NSERC of Canada, the FQRNT of Qu\'ebec and the Universit\'e de Sherbrooke. The second author gratefully acknowledges partial support from CNPq of Brazil. The third author gratefully acknowledges partial support from ANII of Uruguay.}


\begin{abstract}

In this paper, we introduce and study a class of algebras which we call ada algebras. An artin algebra is ada if every indecomposable projective and every indecomposable injective module lies in the union of the left and the right parts of the module category. We describe the Auslander-Reiten components of an ada algebra, showing in particular that its representation theory is entirely contained in that of its left and right supports, which are both tilted algebras. Also, we prove that an ada algebra over an algebraically closed field is simply connected if and only if its first Hochschild cohomology group vanishes.
\end{abstract}

\maketitle

\specialsection*{Introduction}
    \label{Introduction}

Let $A$ be an artin algebra. We are interested in studying the representation theory of $A$, thus the category $\md A$ of finitely generated right $A$-modules. One of the classes of algebras whose representation theory is best understood is that of the quasi-tilted algebras introduced by Happel, Reiten and Smal\o\ in the seminal paper \cite{HRS}. In particular, the ideas and techniques introduced in this paper were used to define and study successfully several generalisations of quasi-tilted algebras, such as shod, weakly shod, laura, left or right supported algebras. For an overview, we refer to the survey \cite{ACLST} or to the more recent \cite{A}.

The objective of present paper is to introduce and study a new class, which we call ada algebras. This also generalises quasi-tilted algebras. Indeed, an artin algebra is quasi-tilted if and only if every indecomposable projective module lies in the so-called left part of the module category, or equivalently if and only if every indecomposable injective module lies in the right part. We say that an algebra is ada if any indecomposable projective and any indecomposable injective lies in the union of these two parts. Ada algebras have the nice property that their representation theory is entirely contained in that of two tilted algebras. Namely, we recall from \cite{ACT,S2} that the left support $A_{\lambda}$ of an artin algebra is the endomorphism ring of the direct sum of all the indecomposable projective modules lying in the left part of $\md A$, and the right support $A_{\rho}$ is defined dually. We prove that the left and right support of an ada algebra are tilted and describe the structure of the module category as in the following theorem.\\
\\
 \textbf{Theorem A} {\it Let $A$ be an ada algebra which is not quasi-tilted. There exists a finite family $(\Gamma_i)_{i=1}^{t}$ of Auslander-Reiten components of $\md A$ which are directed, generalised standard, convex and containing right sections such that:
\begin{itemize}
\item[(a)] $\ind A = \ind A_{\lambda} \cup \ind A_{\rho}$ and each of $A_{\lambda}$ and $A_{\rho}$ is a direct product of tilted algebras.
\item[(b)] If $\Gamma$ is an Auslander-Reiten component of $\md A$ distinct from the $\Gamma_i$, then $\Gamma$ is an Auslander-Reiten component of either $\md A_{\lambda}$ or $\md A_{\rho}$. Moreover
\item[(i)] If $\Hom_{A}(\Gamma, \cup_{i} \Gamma_i ) \neq 0$, then $\Gamma$ is an Auslander-Reiten component of $\md A_{\lambda}$, and,
\item[(ii)] If $\Hom_{A}(\cup{i} \Gamma_i, \Gamma) \neq 0$, then $\Gamma$ is an Auslander-Reiten component of $\md A_{\rho}$.

\end{itemize}}

Furthermore, the portion of the module category of an ada algebra which lies neither in the left nor in the right part is fairly well-understood (see (4.3) below), the structure of the left and right parts being known due to \cite{A}.

 Considering next the case where $A$ is a finite dimensional algebra over an algebraically closed field, we study its simple connectedness. We recall that a triangular algebra $A$ is called simply connected if the fundamental group of any bound quiver presentation of $A$ is trivial, see, for instance \cite{AP}. A well-known problem of Skowro\'nski \cite{S1} links the simple connectedness of $A$ to the vanishing of the first Hochschild cohomology group $HH^{1} (A)$ of $A$ with coefficients in the bimodule ${}_{A}A_{A}$. The equivalence of these conditions holds true for several classes of algebras, and among others for tilted algebras, see \cite{L}. This brings us to our second theorem.\\
\\
\textbf{Theorem B} {\it Let $A$ be an ada algebra over an algebraically closed field. Then $A$ is simply connected if and only if $\HH^{1}(A) = 0$.
Moreover, if this is the case, then the Hochschild cohomology ring $\HH^\bullet(A)$ reduces to the base field}.\\

The paper is organised as follows. After a short preliminary section, we define and study the first properties of ada algebras in section 2. The sections 3 and 4 are occupied with the proof of Theorem A, and section 5 with the proof of Theorem B.


\section{Preliminaries}

\subsection{Notation}
Throughout this paper, all our algebras are basic and connected artin algebras. For an algebra $A$, we denote by $\md A$ its category of finitely generated right modules and by $\ind  A$ a full subcategory of $\md A$ consisting of one representative from each isomorphism class of indecomposable modules. Whenever we speak about a module (or an indecomposable module), we always mean implicitly that it belongs to $\md A$ (or to $\ind  A$, respectively).

Also, all subcategories of $\md A$ are full and so are identified with their object classes. We sometimes consider an algebra $A$ as a category, in which the object class $A_0$ is a complete set $\{e_1, \ldots, e_n \}$ of primitive orthogonal idempotents and the set of morphisms from $e_i$ to $e_j$ is $e_i A e_j$. An algebra $B$ is a \emph{full subcategory} of $A$ if there is an idempotent $e \in A$, sum of some of the distinguished idempotents $e_i$, such that $B = eAe$. It is \emph{convex} in $A$ if, for any sequence $e_i = e_{i_0}, e_{i_1}, \cdots, e_{i_t} = e_j$ of objects in $A$ such that $e_{i_k} A e_{i_{k+1}} \neq 0$ for all $k$, with $0 \leq k < t$, and $e_i, e_j \in B_0$, all $e_{i_{k}}$ lie in $B$. We say that $A$ is \emph{triangular} if there is no sequence $e_i = e_{i_0}, e_{i_1}, \cdots, e_{i_t} = e_i$ of objects in $A$ such that $e_{i_k} A e_{i_{k+1}} \neq 0$ for all $k$, with $0 \leq k < t$. We denote by $P_x$ (or $I_x$, or $S_x$) the indecomposable projective (or injective, or simple, respectively) $A$-module corresponding to the idempotent $e_x$.

Let $\mathcal{C}$ be a subcategory of $\ind  A$. We sometimes write $M \in \mathcal{C}$ to express that $M$ is an object in $\mathcal{C}$. We denote by $\add\mathcal{C}$ the subcategory of $\md A$ with objects the direct sums of summands of modules in $\mathcal{C}$. If $\mathcal{C}, \mathcal{C}'$ are two full subcategories of $\ind  A$, we write $\Hom_{A}(\mathcal{C}, \mathcal{C}')\neq 0$ whenever there exist $M \in \mathcal{C}, M' \in \mathcal{C}'$ such that $\Hom_A(M,M') \neq 0$.

Given a module $M$, we let $\pd M$ (or $\id M$) stand for its projective (or injective, respectively) dimension. The global dimension of $A$ is denoted by $\gd A$.

For an algebra $A$, we denote by $\Gamma(\md A)$ its Auslander-Reiten quiver and $\tau_A = DTr$, $\tau^{-1}_A = TrD$ its Auslander-Reiten translations. For further definitions and facts on $\md A$ or $\Gamma(\md A)$ we refer to \cite{ASS, ARS}.

\subsection{Paths} Let $A$ be an algebra. Given $M,N$ in $ind A$, a \emph{path} from $M$ to $N$ in $\ind  A$ (denoted by $M \leadsto N$) is a sequence of non-zero morphisms\\
$$(\ast) \quad  \quad  M = X_0 \xrightarrow{f_1} X_1 \rightarrow \ldots \xrightarrow{f_t} X_t= N,$$
$(t \geq 1)$ where $X_i \in \ind A$ for all $i$. We then say that $M$ is a \emph{predecessor} of $N$ and $N$ is a \emph{successor} of $M$ (denoted by $M \leq N$).

A path from $M$ to $M$ involving at least one non-isomorphism is a \emph{cycle}. A module $M \in \ind A$ which lies on no cycle is \emph{directed}. If each $f_i$ in $(\ast)$ is irreducible, we say that $(\ast)$ is a \emph{path of irreducible morphisms} or \emph{path in $\Gamma(\md A)$}. A path of irreducible morphisms is \emph{sectional} if $\tau_{A} X_{i+1} \neq X_{i-1}$ for all $i$ with $0 < i < t$.

The left and the right parts of $\md A$ are defined by means of paths. Indeed, the \emph{left part} is the full subcategory of $\ind  A$ with object class
\begin{center}
$\LA = \{ M \in \ind A |$ for any $L$ with $L \leadsto\ M$, we have $\pd L \leq 1 \}$.
\end{center}

Note that $\LA$ is closed under predecessors: if $M \in \LA$ and $L \leadsto M$ then $L \in \LA$. The \emph{right part} $\RA$ is defined dually and is closed under successors.

We need to recall the definitions of Ext-projective and Ext-injective objects. Let $\mathcal{C}$ be a full additive subcategory of $\md A$ closed under extensions (such as $\add \LA$, or $\add \RA$, for instance), then an indecomposable $M \in \mathcal{C}$ is called \emph{Ext-projective} (or \emph{Ext-injective}) in $\mathcal{C}$ if Ext$_{A}^{1}(M,-)|_{\mathcal{C}} = 0$ (or Ext$_{A}^{1}(-,M)|_{\mathcal{C}} = 0$, respectively). It is shown in \cite{AS}(3.4) that $M$ is Ext-injective in $\add \LA$ if and only if $\tau_{A}^{-1} M \notin \LA$ and similarly, $M$ is Ext-projective in $\add \RA$ if and only if $\tau_{A} M \notin \RA$. For further characterisations of these objects, we refer to \cite{ACT}.
\subsection{Left and right section}
A full subquiver $\Sigma$ of a translation quiver $(\Gamma, \tau)$ is called a \emph{right section} if:
\begin{enumerate}
\item $\Sigma$ is acyclic ,
\item for any $x \in \Gamma_0$ such that there exist $y \in \Sigma_{0}$ and a path $y \leadsto x$ in $\Gamma$, there is a unique $n \geq 0$ such that $\tau^{n} x \in \Sigma_{0}$,
\item $\Sigma$ is convex in $\Gamma$.
\end{enumerate}

\emph{Left sections} are defined dually, see \cite{A}. It is shown in \cite{A} that, if $A$ is an artin algebra, and $\Sigma$ is a right section in a generalised standard component of $\Gamma(\md A)$, then $A / \Ann \Sigma$ is a tilted algebra having $\Sigma$ as complete slice \cite{A}(3.6). This notion applies well to the study of the left and right parts. Namely, if $\mathcal{E}$ is the subcategory consisting of all the Ext-projectives in $\add \RA$, and $\Gamma$ is a component of $\Gamma(\md A)$, then:
\begin{itemize}
\item[(a)]If $\Gamma \cap \mathcal{E} = \varnothing$, then either $\Gamma \subseteq \RA$ or $\Gamma \cap \RA = \varnothing$.
\item[(b)] If $\Sigma = \Gamma \cap \mathcal{E} \neq \varnothing$, then $\Sigma$ is a right section of $\Gamma$, convex in $\ind  A$, and moreover $A / \Ann \Sigma$ is a tilted algebra having $\Sigma$ as complete slice, see \cite{A}, Theorem (B).
\end{itemize}
\noindent By component of $\Gamma(\md A)$, we always mean connected component.

\section {Ada algebras: Definition and first properties}

\begin{defn} An artin algebra $A$ is called an \emph{ada algebra} if $A \oplus DA \in \add (\LUR)$.
\end{defn}

Clearly, this is equivalent to requiring that, for every $x \in A_0$, we have both $P_x$ and $I_x$ lying in $\LUR$.

Also, an algebra $A$ is ada if and only if $A^{op}$ is ada. This follows easily from the fact that $D\LA = \mathcal{R}_{A^{op}}$ and $D\RA = \mathcal{L}_{A^{op}}$.

Quasi-tilted algebras are clearly ada. We call \emph{strict} an ada algebra which is not quasi-tilted.

\begin{exs}
\begin{itemize}
  \item [(a)] Let $A$ be a shod algebra \cite{CL1}. Then $\ind  A = \LUR$. Therefore $A$ is ada.
  \item [(b)]Let $A$ be given by the quiver

  $$\xymatrix{%
   \bullet\save[]+<0pt,8pt>*{1}\restore  & \bullet\save[]+<0pt,8pt>*{2}\restore \ar[l]  & \bullet\save[]+<0pt,8pt>*{3}\restore \ar[l]  & \bullet\save[]+<0pt,8pt>*{4}\restore \ar[l] & \bullet\save[]+<0pt,8pt>*{5}\restore \ar[l]}$$

\noindent bound by $\rad^2 A = 0$. Then $P_1, P_2 = I_1, P_3 = I_2$ lie in $\LA$, while $P_4 = I_3$, $P_5 = I_4$ and $I_5$ lie in $\RA$. Then $A$ is a (representation-finite) ada algebra. On the other hand, the one-point extension $A[I_5]$ is not ada.
  \item [(c)] Let $A$ be given by the quiver

$$\xymatrix{%
   \bullet\save[]+<0pt,8pt>*{1}\restore & \bullet\save[]+<0pt,8pt>*{2}\restore \ar@<+2pt>[l]  \ar@<-2pt>[l] & \bullet\save[]+<0pt,8pt>*{3}\restore \ar@<+2pt>[l]  \ar@<-2pt>[l] & \bullet\save[]+<0pt,8pt>*{4}\restore \ar@<+2pt>[l]  \ar@<-2pt>[l]}$$

\noindent bound by $\rad^{2} A =0$. Then $A$ is a (representation-infinite) ada algebra. This example shows that, in contrast to laura algebras \cite{AC1}, an ada algebra may have infinitely many indecomposables which are not in $\LUR$.

\end{itemize}
\end{exs}

Let $P$ denote the direct sum of a complete set of representatives of the isomorphism classes of indecomposable projective $A$-modules lying in $\LA$. Then the algebra $A_{\lambda} = \End P_{A}$ is called the \emph{left support} of $A$, see \cite{ACT, S2}. We recall from \cite{ACT}(2.2) that $A_{\lambda}$ is a full convex subcategory of $A$, closed under successors and that $\LA \subseteq \ind A_{\lambda}$. Moreover, because of \cite{ACT} (2.3), $A_{\lambda}$ (which is not connected in general) is a direct product of quasi-tilted algebras. The \emph{right support} $A_{\rho}$ is defined dually and has dual properties.

\begin{lem}
Let $A$ be an ada algebra, then $A = A_{\lambda} \cup A_{\rho}$.
\end{lem}
\begin{dem}
        Let $x \in A_{0}$. If $P_{x} \in \LA$, then $x \in (A_{\lambda})_{0}$. If not, then $P_x \in \RA$ and the non-zero morphism $P_x \rightarrow I_x$ with image $S_x$ yields $I_x \in \RA$ so that $x \in (A_{\rho})_{0}$.
\end{dem}

\begin{lem}
Let $A$ be an ada algebra, then $A$ is triangular.
\end{lem}

\begin{dem}
Because of \cite{ACT}(2.2)(a), we can write $A$ in triangular matrix form
$A = \left[
       \begin{array}{cc}
         A_\lambda & 0 \\
         M & B \\
       \end{array}
     \right]$.
\noindent Since $A_\lambda$ is a direct product of quasi-tilted algebras, then it is triangular.
On the other hand, let $x \in B_0$, then the indecomposable projective $A$-module $P_x$ does not lie in $\LA$, hence it lies in $\RA$. Now, projectives in $\RA$  are directed because of \cite{A}(6.4). In particular, $B$ is triangular hence so is $A$.
\end{dem}

We have an easy characterisation of ada algebras.

\begin{thm}\label{sub:charac}
An artin algebra $A$ is ada if and only if we have\\ $\ind A = \LA \cup \ind A_{\rho} = \ind A_\lambda \cup \RA$. In particular, if $A$ is ada, then $\ind A = \ind A_\lambda \cup \ind A_{\rho}$.
\end{thm}

\begin{dem}
  Assume first that $A$ is ada, and let $M$ be an indecomposable $A$-module. Suppose that $M \notin \ind A_{\rho}$. Then there exists $x \in A_{0}$ such that $M(x) \neq 0$ and $x \notin (A_{\rho})_{0}$. Thus $I_{x} \notin \RA$ and there exists a non-zero morphism $M \rightarrow I_x$. Since $A$ is ada, then $I_{x} \in \LA$ and so $M \in \LA$. This shows that $\ind A = \LA \cup \ind A_{\rho}$. Similarly, we have $\ind A = \ind A_\lambda \cup \RA$.

  Conversely, assume that these two equalities hold, and let $x \in A_0$, then $P_x \in \RA$ or $P_x \in \ind A_{\lambda}$.
  By definition of $A_{\lambda}$,  this implies $P_x \in \LA$. Therefore $P_x \in \LUR$. Similarly, $I_x \in \LUR$.
\end{dem}

Notice that both conditions $\ind A = \LA \cup \ind A_{\rho}$ and $\ind A = \ind A_{\lambda} \cup \RA$ are necessary for $A$ to be ada.

 We deduce homological properties of ada algebras.

\begin{cor}
Let $A$ be an ada algebra, then
\begin{itemize}
\item [(a)] For any indecomposable module $M$, we have $\pd M \leq 2$ or $\id M \leq 1$.
\item [(b)] $\gd A \leq 4$
\end{itemize}
\end{cor}

\begin{dem}
\noindent (a) This follows from the equality $\ind  A = \ind A_\lambda \cup \RA$ and the fact that $\gd A_{\lambda} \leq 2$ (using that projective $A_{\lambda}$-modules are also projective $A$-modules).\\

\noindent (b) Let $M$ be an indecomposable $A$-module and suppose that $\pd  M \geq 2$. Then there exists a minimal projective resolution
$$0 \rightarrow \Omega^{2}(M) \rightarrow  P_{1} \rightarrow P_{0} \rightarrow M \rightarrow 0$$
and for every indecomposable summand $X$ of $\Omega^{2}(M)$, we have Ext$^{2}_{A}(M,X) \neq 0$. In particular, $\id X \geq 2$. Because of (a), we get $\pd X \leq 2$. This implies that $\pd M \leq 4$.
\end{dem}

\begin{rem}
\begin{itemize}
\item[a)] The bound obtained in (b) above is sharp: indeed, the algebra $A$ of example 2.2(b) has global dimension 4.
\item[b)] Dually, for every $M \in \ind A$, we have $\pd M \leq 1$ or $\id M \leq 2$.
\end{itemize}
\end{rem}

 We now prove that a full subcategory of an ada algebra is ada.

\begin{prop}\label{prop:idem}
Let $A$ be an ada algebra, and $e \in A$ be an idempotent, then $B = eAe$ is ada.
\end{prop}

\begin{dem}
Let $x \in B_0$ and $P_x = e_x B$ denote the corresponding indecomposable projective $B$-module. Then $P_x \otimes_B A \cong e_x A \in \LUR$. Now, because of \cite{AC2}(2.1), we have $\Hom_A (eA,P_x \otimes_{B} A) \in \mathcal{L}_B \cup \mathcal{R}_B$.

But $\Hom_A (eA, P_x \otimes_B A) \cong (P_x \otimes_B A)e \cong e_x A e \cong e_x eAe = e_x B = P_x$.

Then $P_x \in \mathcal{L}_B \cup \mathcal{R}_B$. Similarly, using that $A^{op}$ is ada, we get $I_x \in \mathcal{L}_{B} \cup \mathcal{R}_{B}$.
\end{dem}

 For the notion and main results about split-by-nilpotent extensions, we refer the reader to \cite{AZ}.

\begin{prop}\label{prop:split}
Let $R$ be a split-extension of $A$ by a nilpotent bimodule. If $R$ is ada, then so is $A$.
\end{prop}

\begin{dem}
Let $x \in A_0$, then we clearly have $e_x R_{R} \cong e_x A \otimes_{A} R_{R}$ and $D(Re_{x}) \cong e_x (DR) \cong$,  $\Hom_{A^{op}}(Ae_x, DR) \cong \Hom_{A}(R,D(Ae_x))$. The statement then follows immediately from \cite{AZ}(2.4).
\end{dem}

 Ada algebras also behave well with respect to the skew group algebra construction, see \cite{ARS, ALR}.

\begin{prop}
Let $A$ be an artin algebra, and $G$ be a group acting on $A$ with $|G|$ invertible in $A$. Then the basic algebra $R = A[G]^{b}$ associated to the skew group algebra is ada if and only if $A$ is ada.
\end{prop}

\begin{dem}
Assume first that $A$ is ada, and let $\overline{P}$ be an indecomposable projective $R$-module. Because of \cite{ALR}(4.3), there exists an indecomposable projective summand $P_{A}$ of $\Hom_{R}(R,\overline{P})$ such that $\overline{P}_{R}$ is a direct summand of $P \otimes_{A} R$.

Suppose $P \in \LA$. Because of \cite{ALR}(5.2)(a), we have $P \otimes_{A} R \in \add \mathcal{L}_R$. Therefore $\overline{P} \in \mathcal{L}_R$. Suppose next that $P \in \RA$. Let $X$ be an indecomposable $R$-module such that $\Hom_{R}(\overline{P},X) \neq 0$. We claim that $\id X \leq 1$. Because of \cite{ALR}(4.6), there exist $\sigma \in G$ and an indecomposable summand $M_A$ of $\Hom_{R}(R,X)$ such that $X$ is a summand of ${}^{\sigma}M \otimes_{A} R$ and $\Hom_{A} (P, {}^{\sigma}M) \neq 0$. Because $P \in \RA$, we get $\id {}^{\sigma}M \leq 1$. Since the functor $- \otimes_{A} R : \md A \rightarrow \md R$ is exact and carries injectives to injectives, we get $\id  ({}^{\sigma}M \otimes_{A} R) \leq 1$. Therefore $\id  X \leq 1$, as asserted. Applying \cite{ALR}(1.1) yields $\overline{P} \in \mathcal{R}_{R}$. The proof is entirely similar if we start with an indecomposable injective $R$-module.

Conversely, let $R$ be ada, and $P_{A}$ an indecomposable projective $A$-module. Then there exists an indecomposable projective summand $\overline{P}$ of $P \otimes_{A} R$ such that $P_{A}$ is a direct summand of $\Hom_{R}(R, \overline{P})$.

Suppose $\overline{P} \in \mathcal{L}_R$. Because of \cite{ALR}(5.2)(b), $\Hom_{A}(R, \overline{P}) \in \add \LA$. Therefore $P \in \LA$. Suppose now that $\overline{P} \in \mathcal{R}_{R}$, and let $M$ be an indecomposable $A$-module such that $\Hom_{A} (P,M) \neq 0$. We claim that $\id  M \leq 1$. Because of \cite{RR}, or \cite{ALR}(4.4)(a), we have $\Hom_{R}(\overline{P}, M \otimes_{A}R) \neq 0$. Because of \cite{RR}(1.1 and 1.8), there exists an indecomposable decomposition $M \otimes_{A} R = \oplus_{i=1}^{m} X_i$ such that $\Hom_{R}(R,X_i) = \oplus_{\sigma \in H_i} {}^{\sigma}M$ for some $H_i \subseteq G$. Hence there exists $i$ such that $1 \leq i \leq m$ and $\Hom_{R}(\overline{P},X_i) \neq 0$. Because $\overline{P} \in \mathcal{R}_{R}$, we get $\id   X_i \leq 1$. This implies that, for every $\sigma \in H_i$, we have $\id   {}^{\sigma}M \leq 1$. Therefore $\id   M \leq 1$, as required. Another application of \cite{ALR}(1.1) yields $P \in \RA$. Again the proof is similar if we start with an indecomposable injective $A$-module.
\end{dem}
\section{The module category of an ada algebra}

\subsection{}  Assume $A$ is a strict ada algebra. Then there exists $x \in A_{0}$ such that $P_x \notin \LA$. By definition, $P_x \in \RA$ and is clearly Ext-projective in $\add \RA$. Therefore the set $\Sigma$ of indecomposable Ext-projectives in $\add \RA$ is non-void. Let $\Sigma = \Sigma_{1} \coprod \Sigma_{2} \coprod \cdots \coprod \Sigma_{t}$ where we assume that each $\Sigma_{i}$ is the set of Ext-projectives in $\add \RA$ lying in the same component $\Gamma_{i}$ of $\Gamma (\md A)$. Note that $\Sigma_{i}$ is not necessarily connected.

Because of \cite{A}(6.7), each $\Sigma_{i}$ is a right section in $\Gamma_{i}$, convex in $\ind  A$.
Moreover, $A/\Ann \Sigma_i$ is tilted and has $\Sigma_i$ as a complete slice. The objective of this section is to prove the following theorem.

\begin{thm}\label{theo:components}
Let $A$ be a strict ada algebra. Then
there exists a finite family $(\Gamma_i)_{i=1}^t$ of components of $\Gamma(\md A)$ which are
directed, generalised standard, convex, containing right sections such that, if $\Gamma$ is an Auslander-Reiten component distinct from the $\Gamma_{i}$, then $\Gamma$ is a component of either $\Gamma(\md A_{\lambda})$ or $\Gamma(\md A_{\rho})$ (and, in this latter case, it is contained in $\RA$).
Moreover,
\begin{itemize}
    \item[(i)] if $\Hom_A(\Gamma, \cup_i \Gamma_i) \neq 0$, then
    $\Gamma$ is a component of $\Gamma(\md A_\lambda)$, and
    \item[(ii)] if $\Hom_A(\cup_i \Gamma_i, \Gamma) \neq 0$, then
    $\Gamma$ is a component of $\Gamma(\md A_\rho)$.
    \end{itemize}

\end{thm}
Clearly, the dual statement holds as well: there exists a finite family $(\Gamma'_{j})_{j=1}^{s}$ of directed, generalised standard, convex components of $\Gamma (\md A)$, each containing a left section $\Sigma'_{j}$ consisting of indecomposable Ext-injectives in $\add \LA$, and equipped with the obvious properties. We leave the primal-dual translation to the reader.

We illustrate the theorem with the following example:
\begin{exs} Let $A$ be given by the quiver
 $$\xymatrix{%
    \bullet \ar[r] & \bullet \ar[r] & \bullet  & \bullet\ \ar[l] & \bullet \ar@<+2pt>[l]  \ar@<-2pt>[l] & \bullet\ \ar[l]  }$$
bound by $\rad^{2} A = 0$. The Auslander-Reiten quiver $\Gamma(\md A)$ of $A$ looks as follows.

 $$
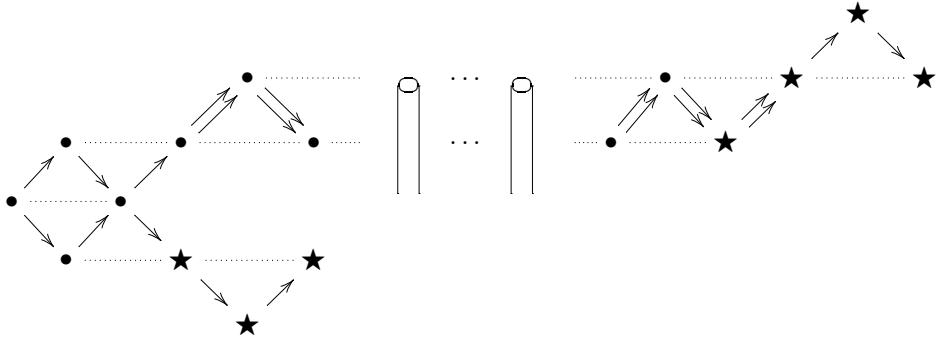
\xymatrix@-1.2pc{%
    &&& &&& &&& &&& &&& \bigstar \ar[rd] & \\
    &&& & \bullet \ar@{.}[rr] \ar@<+2pt>[rd]  \ar@<-2pt>[rd] &&&
      \save[]-<0pt,4pt>*{\put(0,0){\oval(7,5)}}\restore \ar@{-}@<-4pt>[dd] \ar@{-}@<+4pt>[dd] & \ldots &  \save[]-<0pt,4pt>*{\put(0,0){\oval(7,5)}}\restore \ar@{-}@<-4pt>[dd] \ar@{-}@<+4pt>[dd] &  \ar@{.}[rr] & & \bullet \ar@{.}[rr] \ar@<+2pt>[rd]  \ar@<-2pt>[rd] && \bigstar \ar@{.}[rr] \ar[ru] && \bigstar \\
    & \bullet \ar@{.}[rr] \ar[rd] && \bullet \ar@{.}[rr] \ar@<+2pt>[ru]  \ar@<-2pt>[ru] && \bullet \ar@{.}[r] &&& \ldots && \ar@{.}[r] & \bullet \ar@{.}[rr] \ar@<+2pt>[ru]  \ar@<-2pt>[ru] && \bigstar \ar@<+2pt>[ru]  \ar@<-2pt>[ru] &&& \\
    \bullet \ar@{.}[rr] \ar[ru] \ar[rd] && \bullet  \ar[ru] \ar[rd] &&& &&& &&& &&& && \\
    & \bullet \ar@{.}[rr] \ar[ru] && \bigstar \ar@{.}[rr] \ar[rd] && \bigstar &&& &&& &&& && \\
    &&& & \bigstar \ar[ru] &&& &&& &&& &&& \\}$$

\end{exs}
\noindent where we have illustrated the objects of the subcategory $\RA$ by $\bigstar$. Let $\Gamma_1$ denote the postprojective component and $\Gamma_2$ the preinjective component. Then $\Sigma = \Sigma_{1} \cup \Sigma_{2}$ with $\Sigma_1 \subseteq \Gamma_1$ and $\Sigma_2 \subseteq \Gamma_2$. Notice that $\Hom_{A}(\Gamma_1, \Gamma_2) \neq 0$ (and so the components $\Gamma_i$ are not orthogonal). Also, if $\Gamma$ is a regular tube, then $\Hom_{A}(\Gamma_1, \Gamma) \neq 0$ but $\Gamma$ is not contained in $\RA$.

 The proof of Theorem \ref{theo:components} will be split into a series of lemmata.

\begin{lem}\label{lem:successor}
Let $P_x \in \Sigma_i$ be projective. Then every projective successor of $P_x$ lies in the same connected component of $\Sigma_i$.
\end{lem}

\begin{dem}
Assume we have a path $P_x \rightsquigarrow P_y$ with $P_y$ projective. Since $P_x \in \RA$, we have also $P_y \in \RA$. Therefore, $P_y$ is Ext-projective in $\add \RA$ and so there exists $j$ so that $P_y \in \Sigma_j$. By \cite{A}(6.3), the path $P_x \rightsquigarrow P_y$ can be refined to a path of irreducible morphisms and every module on each such refinement is Ext-projective in $\add \RA$. But then, $P_x$ and $P_y$ belong to the same connected component of $\Sigma$. In particular, $i=j$.
\end{dem}

  We denote by $(\Gamma_i)_{\geq \Sigma_i}$ the full subquiver of $\Gamma_i$ consisting of the successors of $\Sigma_i$ (and by $(\Gamma_i)_{\ngeq \Sigma_{i}}$ the full subquiver of $\Gamma_i$ consisting of the non-successors). By definition of $\Sigma$, the successors of $\Sigma_i$ on $\Gamma_i$ are $A_{\rho}$-modules. In fact we have the following result.

\begin{lem}\label{lem:Gammai}
$(\Gamma_i)_{\geq \Sigma_i} = \Gamma_i \cap \RA$
\end{lem}

\begin{dem}
Assume $X \in (\Gamma_i)_{\geq \Sigma_i}$. Then there exist $Y \in \Sigma_i$ and a path $Y \rightsquigarrow X$. Since $Y \in \RA$, we have $X \in \RA$ and so $X \in \Gamma_i \cap \RA$. Conversely, let $X \in \Gamma_i \cap \RA$. Because of \cite{A}(6.6), there exists $m \geq 0$ such that $\tau_{A}^{m} X \in \Sigma_i$. Clearly, $X \in (\Gamma_{i})_{\geq \Sigma_i}$.
\end{dem}

 We have a similar statement for non-successors.

\begin{cor} \label{cor:alambda}
Let $X \in (\Gamma_i)_{\ngeq \Sigma_i}$, then $X \notin \RA$ and $X \in \ind A_{\lambda}$.
\end{cor}

\begin{dem}
The first statement follows from \ref{lem:Gammai}, and the second from \ref{sub:charac}.
\end{dem}

 Since modules in $\Sigma$ are directed (because of \cite{A}(6.4)) we deduce the following statement.

\begin{cor}\label{cor:pred}
Let $X \in \Gamma_i$ be a proper predecessor of $\Sigma$, then $X \notin \RA$ and $X \in \ind A_{\lambda}$.
\end{cor}

\begin{lem}\label{lem:dir}
The modules in $\tau_{A} \Sigma_i$ are directed in $\ind A$.
\end{lem}

\begin{dem}
Since $\Sigma_i$ is acyclic, and $\tau_{A}  \Sigma_i$ contains no injectives, then $\tau_{A} \Sigma_i$ is acyclic. Let $X \in \Sigma_i$ and assume that we have a cycle in $\ind A$

$$\tau_{A} X = M_0 \xrightarrow{f_1} M_1 \rightarrow \ldots \xrightarrow{f_t} M_t=\tau_{A} X.$$

Assume first that none of the $f_j$ factors through an injective module. Then the above cycle induces another one in $\ind A$

$$X = \tau_{A}^{-1} M_0 \rightarrow \tau_{A}^{-1} M_1 \rightarrow \ldots \rightarrow \tau_{A}^{-1} M_t= X.$$

Because of the convexity, this cycle lie inside $\Sigma_i$, thus contradicting the acyclicity of $\Sigma_i$. Therefore, we can assume that there exists $j$ such that $M_j$ is injective. Since, $\tau_{A} X \notin \RA$, we have $M_j \notin \RA$ and thus $M_j \in \LA$. Because of \cite{A}(6.4), $M_j$ is directed, a contradiction.
\end{dem}

\begin{lem}\label{lem:dirc}
For any $i$, $\tau \Sigma_i$ lies in a union of
directed components of $\Gamma(\md A_\lambda)$.
\end{lem}

\begin{dem}
Because of \ref{lem:dir}, $\tau_{A} \Sigma_i$ is directed in $\ind A$, hence it is also directed in $\ind A_\lambda$.

Assume that $X \in \Sigma_i$ is such that $\tau_{A} X$ does not lie in a directed component of $\Gamma(\md A_\lambda)$. Because of the structure of the module category of the quasi-tilted algebra $A_{\lambda}$ (see \cite{C}, \cite{LS}), we have one of two cases:
\begin{enumerate}
  \item $\tau_{A} X$ belongs to an inserted tube or component of type $\mathbb{Z} \mathbb{A}_\infty$ in $\Gamma(\md A_\lambda)$.
  Since $\tau_{A} X$ is directed, there exists a non-directed indecomposable projective $A_\lambda$-module $P$ and a path of irreducible morphisms $\tau_{A} X \rightsquigarrow P$.

Note that $P$ is also projective as an $A$-module and is also not directed in $\ind  A$. In particular, $P \notin \RA$ (by \cite{A}(6.4)). Thus $P \in \LA$ and hence $\tau_{A} X \in \LA$.

 On the other hand, the path $\tau_{A} X \rightsquigarrow P$ of irreducible morphisms contains no injective $A_\lambda$-module, because of the semiregularity of the component. Since any injective $A$-module lying in $\ind A_\lambda$ is also injective as an $A_\lambda$-module, then this path contains no injective $A$-module either. Therefore, we have a path $X \rightsquigarrow \tau^{-1}_{A} P$ of irreducible morphisms. Since $X \in \RA$, then $\tau^{-1}_{A} P \in \RA$. Hence $\tau^{-1}_{A} P \in \Sigma_i$ and $P \in \tau_{A} \Sigma_i$ is directed in $\ind  A$, hence in $\ind  A_{\lambda}$, a contradiction.

  \item $\tau_{A} X$ belongs to a co-inserted tube or component of type $\mathbb{Z} \mathbb{A}_\infty$ in $\Gamma(\md A_\lambda)$. We denote this component by $\Gamma'$.

Recall that $\mathcal{L}_{A_{\lambda}}$ intersects no co-inserted tube or component of type $\mathbb{Z} \mathbb{A}_\infty$. Therefore, no module in $\Gamma'$ belongs to $\mathcal{L}_{A_{\lambda}}$. Because of \ref{sub:charac} and $\LA \subseteq \mathcal{L}_{A_{\lambda}}$, this means that $\Gamma'$ consists entirely of $A_{\rho}$-modules.

We claim that any irreducible morphism $f:Y \rightarrow Z$ between two predecessors of $\tau_{A} X$ in $\Gamma'$ remains irreducible in $\md A$. Indeed, assume that this is not the case, and let $\displaystyle{{g} = \left( \begin{array}{c}
            g_1 \\
            \vdots \\
            g_{t}
            \end{array} \right):Y \rightarrow \oplus_{i=1}^{t} E_i}$ be left minimal
             almost split in $\md A$, where the $E_i$ are assumed indecomposable. Then $f$ factors through $g$, that is, there exists $h = (h_1, \ldots, h_t ): \oplus_{i=1}^{t} E_i \rightarrow Z$ such that $f= \sum_{i=1}^{t} h_i g_i$. Let $i$ be such that $h_i g_i \neq 0$.

Since $Z$ precedes $\tau_{A} X$, then so does $E_i$. Hence $E_i$ is in $\md A_{\lambda}$ by \ref{cor:alambda}. Since so are $Y$ and $Z$, then the left minimal almost split morphism $g$ in $\md A$ remains left minimal almost split in $\md A_{\lambda}$. Consequently, $h$ is a retraction and we are done.

Since $Y,Z$ are predecessors of $\tau_{A} X$ in $\Gamma'$, then they are also indecomposable $A_{\rho}$-modules, and hence $f: Y \rightarrow Z$ remains irreducible in $\md A_{\rho}$.

This implies that the full subquiver $\Gamma'_{\leq \tau_{A} X}$ of all predecessors of $\tau_{A} X$ in $\Gamma'$ is contained in exactly one component $\Gamma$ of $\Gamma (\md A_{\rho})$.

Now, there exist a non-directed injective $A_{\lambda}$-module $I \in \Gamma'$ and a path $I \leadsto \tau_{A} X$ of irreducible morphisms in $\Gamma'$. Because of the previous argument, this path induces a path $I \leadsto \tau_{A} X$ of irreducible morphisms in $\Gamma$. Thus, $\Gamma$ is a component of $\Gamma (\md A_{\rho})$ containing at the same time directed modules (such as $\tau_{A} X$) and non-directed ones (such as $I$) and also a path from a non-directed to a directed module. Using \cite{C}, \cite{LS}, this shows that $\Gamma$ is also a co-inserted tube or component of type
$\mathbb{Z} \mathbb{A}_\infty$ in $\Gamma(\md A_{\rho})$.

Since injective $A_{\rho}$-modules are also injective $A$-modules, there is a non-directed injective $A$-module $J \in \Gamma$ and a path $J \leadsto \tau_{A}X$ in $\ind A_{\rho}$ and therefore in $\ind A$. Since $\tau_{A}X \notin \RA$, then $J \notin \RA$. On the other hand, $J$ is not directed, so $J \notin \LA$, because of \cite{A}(6.4), and this contradicts the hypothesis that $A$ is ada.
\end{enumerate}
\end{dem}

 We may now start the proof of Theorem \ref{theo:components}.

\begin{lem}\label{lem:directed}
Each of the components $\Gamma_i$ is directed and generalised standard and convex in $\ind A$.
\end{lem}

\begin{dem}
Suppose first that we have a cycle in $\ind A$ lying in the component $\Gamma_i$. Since $\Sigma_i$ is a right section, $(\Gamma_{i})_{\geq \Sigma_i}$ is directed, because of \cite{A}(2.2). On the other hand,  $(\Gamma_{i})_{\ngeq \Sigma_i}$ consists of $A_\lambda$-modules, because of \ref{cor:alambda}. We now claim that each connected component $\Gamma$ of $(\Gamma_{i})_{\ngeq \Sigma_i}$ contains at least a module of the form $\tau_{A}X$, with $X \in \Sigma_{i}$.

Assume $\Gamma \cap \tau \Sigma_{i} = \varnothing$. Let $Y \in \Gamma$ (thus, $Y \in \Gamma_{i}$). Since, by definition $\Gamma_{i} \cap \Sigma_{i} \neq \varnothing$ and $\Gamma_{i}$ is connected, then there exists a walk in $\Gamma_{i}$,
$$Y= Y_0 \-- Y_1 \-- \ldots \-- Y_t = X$$
for some $X \in \Sigma_{i}$. We know that $Y$ is not a successor of $\Sigma_{i}$, hence $Y \notin \RA$ while $X \in \RA$. Hence there exists a least $i$ such that $1 \leq i \leq t$ and $Y_0, Y_1, \cdots, Y_{i-1} \notin \RA$ while $Y_i \in \RA$. Then we have an arrow $Y_{i-1} \rightarrow Y_i$. Assume first that $Y_i$ is not projective, then there is an arrow $\tau_{A} Y_i \rightarrow Y_{i-1}$, so $\tau_{A} Y_{i} \notin \RA$. Therefore, $Y_i \in \Sigma_i$. Next, if $Y_i$ is projective, then $Y_{i-1}$ is not injective and so there is an arrow $Y_i \rightarrow \tau^{-1}_{A} Y_{i-1}$. Since $\tau^{-1}_{A} Y_{i-1} \in \RA$ we get $\tau_{A}^{-1} Y_{i-1} \in \Sigma_i$. This establishes our claim. Applying \ref{lem:dirc}, we get that $(\Gamma_i)_{\ngeq \Sigma_i}$ is directed.

This shows that, if we have a cycle in $\Gamma_i$, then it must be of the form

  $$M= M_0 \rightarrow M_1 \rightarrow \ldots \rightarrow M_j \rightarrow \ldots \rightarrow M_t = M$$

\noindent where there exists $j$ such that $M \in (\Gamma_{i})_{\geq \Sigma_i}$ and $M_j \in (\Gamma_{i})_{\ngeq \Sigma_i}$. But now, $M \in (\Gamma_{i})_{\geq \Sigma_i}$ yields $M \in \RA$, and so $M_j \in \RA$, a contradiction to \ref{lem:Gammai}. This shows that $\Gamma_i$ is directed.

  Now, we assume that $\Gamma_i$ is not generalised standard and let $L$, $M \in \Gamma_i$ be such that $\rad_A^\infty(L, M) \neq 0$. Since $(\Gamma_{i})_{\geq \Sigma_i}$ is generalised standard, because of \cite{A}(3.2), and $(\Gamma_{i})_{\ngeq \Sigma_i}$ also, because it is part of a directed, hence generalised standard component of the Auslander-Reiten quiver of the quasi-tilted algebra $A_{\lambda}$, then we must have $L \in (\Gamma_{i})_{\ngeq \Sigma_i}$ and $M \in (\Gamma_{i})_{\geq \Sigma_i}$. Let $f \in \rad_A^\infty(L, M)$ be non-zero. For any $t\geq 0$, the morphism $f$ induces a path in $\ind A$

  $$L \xrightarrow{g_t} M_t \xrightarrow{f_t} \ldots \rightarrow M_1 \xrightarrow{f_1} M_0=M$$

\noindent with $f_1, \ldots, f_t$ irreducible, $g_t \in \rad_A^\infty(L, M_t)$ and $f_1 \ldots f_t g_t \neq 0$. Therefore, there exists $t$ such that $M_t \in (\Gamma_{i})_{\ngeq \Sigma_i}$ and $\rad_A^\infty(L, M_t) \neq 0$, a contradiction to the fact that $(\Gamma_{i})_{\ngeq \Sigma_i}$ is generalised standard.

It remains to prove the convexity of $\Gamma_i$. Assume that we have a path in $\ind A$:

  $$M = M_0 \xrightarrow{f_1} M_1 \rightarrow \ldots \xrightarrow{f_t} M_t=N$$

  \noindent with $M$, $N \in \Gamma_i$ and $M_1, \ldots, M_{t-1} \notin \Gamma_i$ (thus $t \geq 2$).
  Then, $f_t \in \rad _A^\infty (M_{t-1}, N)$. Suppose first that $N \in (\Gamma_{i})_{\geq \Sigma_i}$ then, for any $s \geq 0$, we have a path in $\ind A$

 $$M_{t-1} \xrightarrow{h_s} N_s \xrightarrow{g_s} \ldots \rightarrow N_1 \xrightarrow{g_1} N_0=N$$

\noindent with $g_1, \ldots, g_s$ irreducible and $h_s \in \rad_A^\infty(M_{t-1}, N_s)$ such that $h_s g_s \ldots g_1 \neq 0$. Then there exists $s$ such that $N_s \in (\Gamma_{i})_{\ngeq \Sigma_i}$.

We may thus suppose from the start that $N \in (\Gamma_{i})_{\ngeq \Sigma_i}$. In particular, $N \notin \RA$ and thus $M \notin \RA$ and they are $A_\lambda$-modules because of \ref{cor:alambda}. We claim that all $M_j$ are $A_\lambda$-modules. Indeed, if this is not the case, by \ref{sub:charac} there exists $M_j \in \RA$, a contradiction. Then the given path consists entirely of $A_{\lambda}$-modules, with $M,N \in (\Gamma_i)_{\ngeq \Sigma_i}$.
The conclusion then follows from the fact that $(\Gamma_i )_{\ngeq \Sigma_i}$ is part of a directed component, hence convex component of $\Gamma(\md A_{\lambda})$.
\end{dem}

 Recall that an artin algebra $A$ is \emph{laura} if the class $\ind A \setminus (\LUR)$ contains only finitely many objects \cite{AC1}. A laura algebra which is not quasi-tilted always has a unique Auslander-Reiten component which is non-semiregular and faithful. The algebra $A$ is called \emph{weakly shod} \cite{CL2} if this component is directed.

\begin{cor}
Let $A$ be a strict ada algebra. If $A$ is laura, then it is weakly shod.
\end{cor}

\begin{dem}
Let $\Gamma$ be the faithful non-semiregular component of $\Gamma(\md A)$. Since $A$ is strict, there exists a projective $A$-module $P_x$ such that $P_x\in \RA\setminus\LA$. Because $\Gamma$ is faithful, there exists $M\in \Gamma$ such that $\Hom_A(P_x,M)\neq 0$ and so $M\in\RA\setminus\LA$. This shows that $\Gamma\cap\RA\neq 0$ and that $\Gamma \nsubseteq \LA$. Dually $\Gamma \nsubseteq \RA$.

Because of \cite{A}, Theorem B, the intersection of $\Gamma$ with the class $\Sigma$ of indecomposable Ext-projectives in $\add \RA$ is a right section of $\Gamma$. Since $\Gamma = \Gamma_i$ is directed because of \ref{lem:directed}, we get that $A$ is weakly shod.
\end{dem}

 The proof of Theorem \ref{theo:components} will be completed once we prove the following lemma
\begin{lem}
Let $A$ be a strict ada algebra. If $\Gamma$ is a component of $\Gamma(\md A)$ distinct from the $\Gamma_i$, then $\Gamma$ is a component of either $\Gamma(\md A_{\lambda})$ or $\Gamma(\md A_{\rho})$ (and in this latter case, it is contained in $\RA$). Moreover, we have either

    \begin{itemize}
    \item[i)] If $\Hom_A(\Gamma, \cup_i \Gamma_i) \neq 0$ then
    $\Gamma$ is a component of $\Gamma(\md A_\lambda)$, or
    \item[ii)] If $\Hom_A(\cup_i \Gamma_i, \Gamma) \neq 0$ then
    $\Gamma$ is a component of $\Gamma(\md A_\rho)$
    \end{itemize}
\end{lem}

\begin{dem}
Because $\Gamma \neq \Gamma_i$ for all $i$, we have $\Gamma \cap \Sigma = \emptyset$. Because of \cite{A}(Theorem B), we get that either $\Gamma \subseteq \RA$ or $\Gamma \cap \RA = \emptyset$. In the first case, clearly, $\Gamma$ is a component of $\Gamma (\md A_{\rho})$ contained in $\RA$. We claim that, if $\Gamma \cap \RA = \emptyset$, then $\Gamma$ is a component of $\Gamma (\md A_{\lambda})$. It suffices to prove that each $X \in \Gamma$ is an $A_{\lambda}$-module. Now, if this is not the case, then there exists an indecomposable projective $P \notin \LA$ such that $\Hom_A(P,X) \neq 0$. But then $P \in \RA$ and so $X \in \RA$, a contradiction which establishes our claim.

Now, assume that $\Hom_A (\Gamma, \cup _i \Gamma_i) \neq 0$ and $\Gamma$ is not a component of $\Gamma (\md A_{\lambda})$. Let $X \in \Gamma$ be not an $A_{\lambda}$-module. Then there exists an indecomposable projective $A$-module $P \notin \LA$ such that $\Hom_A (P,X) \neq 0$. As above, $X \in \RA$ and so $\Gamma \cap \RA \neq \emptyset$. Because of  \cite{A}(Theorem B), we have $\Gamma \subseteq \RA$.

Since $\Hom_A(\Gamma, \cup_i \Gamma_i) \neq 0$, there exist $M \in \Gamma$ and $N \in \Gamma_i$ for some $i$ such that $\Hom _A(M, N) \neq 0$.
Since $M \in \RA$, thus $N \in \RA$. Because of \ref{lem:Gammai}, we have $N \in (\Gamma_i)_{\geq \Sigma_i}$. Since $\Gamma \neq \Gamma_i$, we have $\Hom _A(M, N)= \rad ^\infty_A(M, N) \neq 0$. Thus, for any $s \geq 0$, there exists a path in $\ind A$

  $$M \xrightarrow{h_s} N_s \xrightarrow{g_s} \ldots \rightarrow N_1 \xrightarrow{g_1} N_0=N$$

\noindent with $g_1, \ldots, g_s$ irreducible and $h_s \in \rad_A^\infty(M, N_s)$ such that $g_1 \ldots  g_s h_s\neq 0$. Therefore, there exists $s$ such that $N_s \in (\Gamma_{i})_{\ngeq \Sigma_i}$. But then $N_s \in \RA$, a contradiction to \ref{lem:Gammai}. This completes the proof of i).

Finally, assume similarly that $\Hom_A (\cup_{i} \Gamma, \Gamma) \neq 0$ and $\Gamma$ is not a component of $\Gamma (\md A_{\rho})$. In particular, $\Gamma$ is not contained in $\RA$ and since moreover $\Gamma \cap \Sigma = \emptyset$, we deduce from \cite{A}, Theorem B, that $\Gamma \cap \RA = \emptyset$.

By hypothesis, there exist $i$, $M \in \Gamma_i$ and $X \in \Gamma$ such that $\Hom_A (M,X) \neq 0$. If $M \in (\Gamma_i)_{\geq \Sigma _i}$, then $M \in \RA$ by \ref{lem:Gammai}, so that $X \in \RA$, a contradiction. Therefore, $M$ is not a successor of $\Sigma_i$. We then consider two cases.

Suppose first that $( \Gamma_i)_{\ngeq \Sigma_i}$ contains no injective. In this case, $\Sigma_i$ is a section in the directed component $\Gamma_i$, because of \cite{A}(2.3) and moreover $\Gamma_i$ is the connecting component of the tilted algebra $A / \Ann \Sigma_i$, and $\Sigma_i$ is a complete slice, because of \cite{A}(3.6). Now, observe that $\Sigma_i \subseteq \RA$, so $(\Gamma_i)_{\geq \Sigma_i} \subseteq \RA$, thus $(\Gamma_i)_{\geq \Sigma_i}$ consists of $A_{\rho}$-modules. Since $\Sigma_i$ cogenerates $(\Gamma_i)_{\ngeq \Gamma_i}$, then $(\Gamma_i)_{\ngeq \Sigma_i}$ also consists of $A_{\rho}$-modules. In particular, $A /\Ann \Sigma_i$ is a connected component of $A_{\rho}$. Because $\Sigma_i$ is a complete slice, $M \in \Sigma_i$ is not a successor of $\Sigma_i$ if and only if $M$ is a predecessor of $\Sigma_i$. Therefore $\rad^{\infty}_{A}(M,X) \neq 0$ gives, for any $t \geq 0$, a path in $\ind  A$
$$M = M_0 \xrightarrow{f_1} M_1 \rightarrow \ldots \xrightarrow{f_t} M_t  \xrightarrow{g_t} X$$
\noindent where the $f_i$ are irreducible and $g_t \in \rad^{\infty}_{A}(M_t,X)$ is such that $g_t f_t \cdots f_1 \neq 0$. Let $t \geq 0$ be such that $M_t$ is a successor of $\Sigma_i$, then $M_t \in \RA$, hence $X \in \RA$ and we get a contradiction in this case.

Suppose next that $( \Gamma_i)_{\ngeq \Sigma}$, contains an injective $A$-module $I$. Because of \ref{cor:alambda}, we have $I \notin \RA$. Hence $I \in \LA$ and so is Ext-injective in $\add \LA$. Using the notation in \ref{theo:components}, this shows that the Ext-injectives in $\add \LA$ form a left section $\Sigma'_{j}$ in some component $\Gamma'_{j}$. Note that $\Gamma'_{j} = \Gamma_i$. Since $\rad^{\infty}_{A}(M,X) \neq 0$, there exists, for each $t \geq 0$, a path in $\ind  A$
$$M = M_0 \xrightarrow{f_1} M_1 \rightarrow \ldots \xrightarrow{f_t} M_t  \xrightarrow{g_t} X$$
\noindent where the $f_i$ are irreducible and $g_t \in \rad^{\infty}_{A}(M_t,X)$ is such that $g_t f_t \cdots f_1 \neq 0$. Let $t \geq 0$ be such that $M_t$ is a proper successor of $\Sigma'_{j}$. Because of
\ref{lem:Gammai}, this gives $M_t \notin \LA$. Therefore, $X \notin \LA$. This shows that $\Gamma$ contains at least an indecomposable $X$ which is not in $\LA$. Now, we claim that $\Gamma \cap \LA = \emptyset$. By induction, it suffices to show that no neighbour $Y$ of $X$ belongs to $\LA$. If there is an arrow $X \to Y$, then $X \notin \LA$ implies $Y \notin \LA$. Assume that we have an arrow $Y \to X$ and that $Y \in \LA$. We claim that in this case $Y$ is Ext-injective in $\add \LA$. This is obvious if $Y$ is injective, and, if it is not, then there is an arrow $X \to \tau^{-1}_{A} Y$ so that $\tau^{-1}_{A} Y \notin \LA$ and again $Y$ is Ext-injective in $\add \LA$. In particular, $\Gamma = \Gamma'_{l}$ for some $l$ and $Y \in \Sigma'_{l}$. Now there exists a non-zero morphism $g_s \in \rad^{\infty}_{A}(M_s,X)$.
This morphism factors through $\Sigma'_{l}$ (because $X$ is a successor of $\Sigma'_{l}$). Then $\Sigma'_{l} \subseteq \LA$ yields $M_s \in \LA$ and this is a contradiction. Therefore $Y \notin \LA$. This shows that $\Gamma \cap \LA = \emptyset$. Because of \ref{sub:charac}, $\Gamma$ consists of $A_{\rho}$-modules and hence is a component of $\Gamma (\md A_{\rho})$.
\end{dem}

\section{The supports of an ada algebra}

Throughout this section, we let $A$ be a strict ada algebra.

\begin{prop}\label{lem:tilted}
Each of $A_\lambda$ and $A_\rho$ is a direct product of tilted algebras.
\end{prop}

\begin{dem}
Indeed, assume that $B$ is a connected component of $A_\lambda$ and is not tilted. Since $A$ is strict, we have $B \neq A$ and so there exist an indecomposable $B$-module $X$ and an irreducible morphism $X \to P_x$ with $P_x$ an indecomposable projective $A$-module which is not a $B$-module. Since $X$ is isomorphic to an indecomposable summand of $\rad_A(P_x)$, then $P_x \notin \LA$ hence $P_x \in \RA$ and therefore is Ext-projective in $\add \RA$.

We claim that $X$ is a directed $A$-module. Indeed, $X$ is not injective, so we have an arrow $P_x \to \tau^{-1}_A X$ and then we have two cases. If $X \notin \RA$ then $\tau_{A}^{-1}X \in \RA$ yields $\tau_{A}^{-1} X \in \Sigma$ and so $X \in \tau_{A} \Sigma$ is a directed $A$-module. If $X \in \RA$, then $X \in \Sigma$ and so is again directed. In fact, it follows from \ref{lem:directed} that $X$ lies in a directed component of $\Gamma(\md A)$ and \ref{lem:dirc} that it lies in a directed component of $\Gamma(\md B)$. Since $B$ is quasi-tilted but not tilted, then this is the postprojective or the preinjective component of $\Gamma(\md B)$.

Let $e = e_x + \sum_{y \in B_0} e_y$. Then $A'= eAe$ is ada, because of \ref{prop:idem} and is a one-point extension of $B$. Because of \ref{prop:split}, we may assume that $A'= B[X]$.

Assume first that $X$ lies in the postprojective component of $\Gamma(\md B)$. Let $P'_x$ be the indecomposable projective $A'$-module corresponding to the point $x$. Then, considering $P'_x$ as an $A$-module under the standard embedding of $\md A'$ into $\md A$, we have an epimorphism $P_x \to P'_x$. Since $P_x \in \RA \setminus \LA$, then $P'_x \in \RA \setminus \LA$ as well. Applying \cite{AC2}(2.1), we get $P'_x \in \mathcal{R}_{A'}$. On the other hand, since $B$ is quasi-tilted but not tilted, there exists a non-directed indecomposable projective $B$-module $P_y$ lying in an inserted tube or component of type $\mathbb{ZA}_\infty$. Note that $y$ is a source in $B$ and hence also is $A'$. Thus $P_y = P_y'$ is a non-directed indecomposable projective $A'$-module. On the other hand, $P'_x$ lies in the postprojective component of $\Gamma(\md A')$. We claim that there exists a path $P'_x \leadsto P'_y$ in $\md A'$. Indeed, since $B$ is connected and $y$ is a source, there exists $z \in B_0$ such that $P'_z$ lies in the postprojective component of $\Gamma(\md A')$ and a non-zero morphism $f : P'_z \to P'_y$. Since $f \in \rad_{A'}^{\infty}(P'_z,P'_y)$, there exists, for any $t \geq 0$, a path in $\ind A$
$$P'_z = M_0 \xrightarrow{f_1} M_1 \rightarrow \ldots \xrightarrow{f_t} M_t  \xrightarrow{g_t} P'_y$$
with the $f_i$ irreducible and $g_t \in \rad_{A'}^{\infty}$ ($M_t, P'_y$) such that $g_t f_t \ldots f_1 \neq 0$.

Let $t$ be such that $M_t$ is a successor of $P'_x$. This yields the required path $P'_x \leadsto P'_y$ in $\md A'$. But we have already seen that $P'_x \in \mathcal{R}_{A'}$, a contradiction because $P'_y$ is not directed.

Therefore, we may assume $X$ to lie in the preinjective component of $\Gamma(\md B)$. Now, since $B$ is quasi-tilted but not tilted, there exists a non-directed indecomposable injective $B$-module $I_y$ lying in a co-inserted tube or component of type $\mathbb{ZA}_\infty$. Because $A'= B[X]$ and $X$ is preinjective, then $I_y$ is also an injective $A'$-module. However, we have $P'_x \in \mathcal{R}_{A'}$, and there exists a non-sectional path $I_y \leadsto X \rightarrow P'_x$. Because of \cite{AC2}(1.5), this implies that $I_y \notin \mathcal{R}_{A'}$. The algebra $A'$ being ada, we get $I_y \in \mathcal{L}_{A'}$ a contradiction, because $I_y$ is not directed. The proof is now complete.
\end{dem}

 It follows from \ref{theo:components} and \ref{lem:tilted} that, if $A$ is an ada algebra, then we have a good description of the indecomposable modules (or components) lying in $\LUR$: these are modules (or components) over one of the tilted algebras $A_{\lambda}$ and $A_{\rho}$. We now wish to describe those modules which do not belong to $\LUR$. As in \ref{theo:components}, we denote by $\Sigma$ the class of Ext-projectives in $\add \RA$ and by $\Sigma'$ the class of Ext-injectives in $\add \LA$.

\begin{lem}\label{lem:nsectional}
Let $A$ be a strict ada algebra and $X$ an indecomposable $A$-module not lying in $\LUR$. Then there exist an indecomposable projective module $P \in \Sigma$ and a non-sectional path $X \leadsto P$.
\end{lem}
\begin{dem}
Indeed, since $X \notin \RA$, then there exists a path $X \leadsto Y$ in $\ind A$ where $Y$ is such that $\id Y > 1$. Hence there exists an indecomposable projective $A$-module $P$ such that we have a path $X \leadsto Y \rightarrow \ast \rightarrow \tau^{-1}_{A} Y \rightarrow P$ in $\ind A$. Since $X \notin \LA$, we also have $P \notin \LA$. Therefore $P \in \RA$ and so $P \in \Sigma$.
\end{dem}

 Now, notice that $C = A_{\lambda}\cap A_{\rho}$ is a full convex subcategory of $A_{\lambda}$ (or $A_{\rho}$) and therefore is tilted, because of \cite{H2}(III.6.5).

\begin{prop}
Let $A$ be a strict ada algebra, and $X$ be an indecomposable $A$-module. The following conditions are equivalent.
\begin{itemize}
\item[(a)] $X \notin \LUR$.
\item[(b)] There exist $P \in \Sigma$ projective, $I \in \Sigma'$ injective and two non-sectional paths $I \leadsto X$ and $X \leadsto P$.
\item[(c)] $X$ is a proper predecessor of $\Sigma$ and a proper successor of $\Sigma'$.
\end{itemize}

Moreover, if this is the case, then $X$ is an indecomposable $C$-module, generated by $\Sigma'$ and cogenerated by $\Sigma$.
\end{prop}
\begin{dem}
That (a) implies (b) follows from \ref{lem:nsectional} and its dual.
That (b) implies (c) follows from \cite{A}(6.3), because the given paths are non-sectional. Finally, assume that (c) holds. Since $X$ is a proper predecessor of $\Sigma$, then there exists a non-sectional path from $X$ to some $M \in \Sigma$. Because of \cite{A}(6.3), this implies that $X \notin \RA$. Similarly, $X \notin \LA$.

Now, if this is the case, then $X$ being a proper predecessor of $\Sigma$ implies $X \in \ind A_{\lambda}$, because of \ref{cor:pred}. Similarly, $X \in \ind A_{\rho}$. Therefore $X \in \ind C$. The statements about generation and cogeneration follow from the fact that there exist neither projectives nor injectives lying strictly between $\Sigma'$ and $\Sigma$.
\end{dem}


\section{Hochschild cohomology and simple connectedness}

Throughout this last section, all our algebras are finite dimensional algebras over an algebraically closed field $k$.

Let $A$ be ada. We recall from \cite{AL} that an indecomposable projective $P_x \in \RA$ is called a \emph{maximal projective} if it has no projective successor. We then say that $A$ is a \emph{maximal extension} of $B=A \setminus \{x\}$. Denoting by $M$ the radical of $P_x$, we have $A=B[M]$. We shall prove in \ref{prop:maxfilt} below that any strict ada algebra may be written as a maximal extension of another ada algebra.

\begin{lem}\label{lem:maxex}
Let $A=B[M]$ be a maximal extension. Then for every $i \geq 1$, we have $\Ext_{B}^{i}(M,M) = 0$.
\end{lem}
\begin{dem}
Same as \cite{AL}(2.3).
\end{dem}

Let $\HH^{i}(A)$ denote the $i^{th}$ Hochschild cohomology group of $A$ with coefficients in the bimodule ${}_{A}A_{A}$ (see \cite{H1} for details). It is shown in \cite{H1}(5.3) that, if $A = B[M]$, then there exists a long exact sequence
$$0 \rightarrow \HH^{0}(A) \rightarrow \HH^{0}(B) \rightarrow \End M/k \rightarrow \HH^{1}(A) \rightarrow \HH^{1}(B) \rightarrow \Ext^{1}_{B}(M,M) \rightarrow \cdots$$
$$\cdots \rightarrow \HH^{i}(A) \rightarrow \HH^{i}(B) \rightarrow \text{Ext}^{i}_{B}(M,M) \rightarrow \cdots$$
We refer to this sequence in the sequel as \emph{Happel's sequence}. We also recall that the extension point $x$ is called \emph{separating} if the number of indecomposable summands of $\rad P_x$ equals the number of connected components of $B = A\setminus \{x\}$, see, for instance \cite{AP}.

\begin{lem}\label{lem:maxexada}
Let $A=B[M]$ be an ada maximal extension. Then:
\begin{itemize}
\item [(a)] There exists an exact sequence $$0 \rightarrow \HH^{0}(A) \rightarrow \HH^{0}(B) \rightarrow \End M/k \rightarrow \HH^{1}(A) \rightarrow \HH^{1}(B) \rightarrow 0$$
\item [(b)] For any $i \geq 2$, we have $\HH^{i} (A) \cong \HH^{i} (B)$.
\item [(c)] $\HH^{1} (A) \cong \HH^{1} (B)$ if and only if the extension point is separating.
\end{itemize}
\end{lem}

\begin{dem}

The statements (a) and (b) follow from Lemma \ref{lem:maxex} and Happel's sequence. We proceed to prove (c). The surjective morphism $\HH^{1}(A) \rightarrow \HH^{1}(B)$ has kernel with dimension equal to $$\dim_{k}(\text{End} M/k) - \dim_{k} \HH^{0}(B) + \dim_{k} \HH^{0}(A) = \dim_{k} \End M - \dim_{k} \HH^{0}(B)$$ because $A$ is connected. Therefore, $\HH^{1}(A) \cong \HH^{1}(B)$ if and only if $\dim_{k} \End M$ equals the number of connected components of $B$, and this is the case if and only if the extension point $x$ is separating and $M$ is a direct sum of bricks. Because of Theorem \ref{theo:components}, every indecomposable projective lying in $\RA$ belongs to a directed generalised standard component. Therefore, every indecomposable summand of $M$ is a brick. The statement follows.
\end{dem}

\begin{rem}
In particular, we proved that the module $M$ is separated, see \cite{AP} for the definition.
\end{rem}

 A triangular algebra $A$ is called \emph{simply connected} if, for every presentation $A \cong kQ/I$ of $A$ as a bound quiver algebra, the fundamental group of $(Q,I)$ is trivial, see \cite{S1, AP}. Let $A=B[M]$ where we denote by $x$ the extension point. We fix a presentation of $A$ and consider the induced presentation of $B$. Let $\sim$ be the least equivalence relation on the arrows of source $x$ such that $\alpha_1 \sim \alpha_2$ if there exists a minimal relation of the form $\lambda_1 \alpha_1 v_1 + \lambda_2 \alpha_2v_2 + \sum_{j \geq 3} \lambda_j w_j$. Let $t$ be the number of equivalence classes of arrows of source $x$ under this relation. For each $i$, with $1 \leq i \leq t$, let $l(i)$ be the number of tuples of paths $(u_1, v_1, \ldots, u_n, v_n)$ such that there are minimal relations of the forms $\lambda_{1,1}\alpha_1 u_1 + \lambda_{2,1} \alpha_n v_n + \sum_{j \geq 3} \lambda_{j,1} w_{j,1}, \lambda_{1,2} \alpha_1 v_1 + \lambda_{2,2}\alpha_{2} u_2 + \sum_{j \geq 3} \lambda_{j,2} w_{j,2}, \cdots$ where $\alpha_1, \cdots, \alpha_n$ are distinct arrows in the same equivalence class, see \cite{AP}(2.4).

\begin{lem}\label{lem:sc}
Let $A$ be a strict ada algebra.
\begin{itemize}
\item [(a)] If $B$ is a direct product of simply connected algebras, then $A$ is simply connected if and only if the extension point is separating.
\item [(b)] If $A$ is a simply connected strict ada maximal extension, then $B$ is a direct product of simply connected algebras.
\end{itemize}
\end{lem}

\begin{dem}
\noindent (a) This statement follows from \cite{ABCN}(3.6).\\
\noindent (b) Let $B \cong kQ_B/I'$ be an arbitrary presentation of $B$, then there exist a presentation $A \cong kQ_A/I$ of $A$ such that $I \cap kQ_B = I'$. Because of \cite{AP}(2.4) it suffices to show that $l(i) =0$ for all $i$. However, if $l(i) \neq 0$ for some $i$, then there exists a tuple of paths $(u_1, v_1, \cdots, u_n, v_n)$ and a full subcategory $C$ of $A$ which is a split extension of a subcategory $D$ of the form

$$
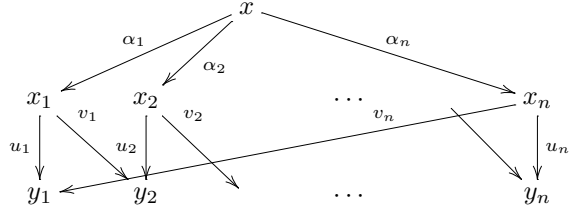
\xymatrix{ & & x \ar[lld]_{\alpha_1} \ar[ld]^{\alpha_2} \ar[rrrd]^{\alpha_n}
& & & \\ x_1 \ar[d]_{u_1} \ar[rd]^(.3){v_1} & x_2 \ar[d]_{u_2}
\ar[rd]^(.3){v_2} & & \cdots & \ar[rd] & x_n \ar[d]^ {u_n}
\ar[llllld]_(.3){v_n} \\ y_1 & y_2 & & \cdots & & y_n}$$\\

\noindent (indeed, there might be in $C$ additional arrows from some $y_i$ to some $y_j$). We denote respectively by $P_x$, $P'_x$, $P''_x$ the indecomposable projective module corresponding to $x$ in $\md A$, $\md C$ and $\md D$.
Then $P'_x = P''_x \otimes_D C$ and we have an epimorphism from $P_x$ to $\overline{P}'_x$
where $\overline{P}'_x = P'_x \otimes_{C} A$. Now, $P_x \in \RA \setminus \LA$ (because $A$ is strict), hence $\overline{P}'_x \in \RA \setminus \LA$. But then, because of \cite{AC2}(2.1), $P'_x \in \mathcal{R}_{C}$.
Hence, because of \cite{AZ}(2.4), we have $P''_x \in \mathcal{R}_D$. However, $\rad P''_x$ is a simple homogeneous module over the hereditary full subcategory of $D$ with class of objects $D \setminus \{x\}$. In particular, $\rad P''_x$ is not directed in $\ind  D$, hence neither is $P''_x$. This however contradicts the fact that $P''_x \in \mathcal{R}_{D}$ (and \cite{A} (6.4)). Therefore $l(i) = 0$ for all $i$ as asserted and so $B$ is a direct product of simply connected algebras.
\end{dem}

 We say that an ada algebra is of \emph{tree type} if the orbit graph (see, for instance, \cite{BG} or \cite{AL}(4.1)) of each of the $\Gamma_i$ is a tree.

\begin{lem}\label{lem:ext}
Let $A = B[M]$ be an ada maximal extension. Then $A$ is of tree type if and only if $B$ is of tree type and the extension point is separating.
\end{lem}

\begin{dem}
Same as \cite{AL}(4.1).
\end{dem}

 A sequence of ada algebras of the form
$$A_{\lambda} = A_0 \varsubsetneqq A_1 \varsubsetneqq \cdots \varsubsetneqq A_m =A$$
is called a \emph{maximal filtration} of $A$ provided that for each $i$, with $1 \leq i \leq m$, there exists an $A_{i-1}$-module $M_i$ such that $A_i = A_{i-1}[M_i]$ is a maximal extension.

\begin{prop}\label{prop:maxfilt}
Let $A$ be a strict ada algebra. Then $A$ admits a maximal filtration.

\end{prop}

\begin{dem}
Since $A$ is strict, there exists an indecomposable projective in $\RA$ which is not in $\LA$. Since every such projective is directed, because of \cite{A}(6.4), there exists (at least) a maximal projective $P_x$. Let $A = B[M]$ where $B = A \setminus \{x\}$ and $M = \rad P_x$. Because of  \ref{prop:idem}, $B$ is also an ada algebra. If $B$ is not strict, then every indecomposable projective $B$-module lies in $\LA \cap \ind B = \mathcal{L}_B \subseteq \LA$ and so $B = A_{\lambda}$. Otherwise, we apply induction.
\end{dem}

\begin{cor}\label{cor:coho}

Let $A$ be a strict ada algebra, then
\begin{itemize}
\item[(a)] $\HH^{1}(A) = 0$ if and only $\HH^{1}(A_{\lambda}) = 0$ and each of the extension points of a maximal filtration is separating.
\item[(b)] $\HH^{i}(A) = 0$ for all $i \geq 2$.
\end{itemize}
\end{cor}

\begin{dem}
\noindent (a) This follows immediately from \ref{lem:ext} and \ref{lem:maxexada}.\\
\noindent (b) Follows from \ref{lem:ext} and \ref{lem:maxexada}, using that $A_{\lambda}$ is tilted and \cite{H3}, Theorem 2.2.
\end{dem}

We also have the immediate corollary.
\begin{cor}\label{cor:tree}
Let $A$ be a strict ada algebra.
Then $A$ is of tree type if and only if $A_{\lambda}$ is of tree type and each of the extension points in a maximal filtration is separating.\hfill $\Box$
\end{cor}

 We are now in a position to prove our main result of this section.
\begin{thm}\label{theo:main}
Let $A$ be an ada algebra. The following are equivalent:
\begin{itemize}
\item [(a)] $A$ is simply connected.
\item [(b)] $\HH^{1}(A) = 0$
\item [(c)] $A$ is of tree type.
\end{itemize}
\end{thm}

\begin{dem}
We may assume that $A$ is strict ada.

Assume first that $\HH^{1}(A) = 0$. Because of \ref{cor:coho}(a), we have $\HH^{1}(A_{\lambda})=0$ and each of the extension points in a maximal filtration is separating. Because of \cite{L}, $\HH^{1}(A_{\lambda})=0$ if and only if $A_{\lambda}$ is a direct product of simply connected algebras. Applying \ref{lem:sc}(a) and induction, we get that $A$ is simply connected.

Conversely, assume that $A$ is a simply connected ada algebra. Therefore there exists a maximal projective $P_x \in \RA$, such that $A=B[M]$ is a maximal extension where, as usual, $B=A \setminus \{x\}$ and $M= \rad P_x$. Now, $x$ is a source in $A$, hence, by \cite{AP}(2.6), $x$ is separating. On the other hand, because of \ref{lem:sc}(b), $B$ is a direct product of simply connected algebras. Hence, inductively, $\HH^{1}(B) = 0$. Applying \ref{lem:maxexada}(c), we get $\HH^{1}(A) = 0$.

The equivalence with condition (c) is proved in the same way using \ref{cor:tree}, and the fact proved in \cite{L}, that $A_{\lambda}$ is of tree type if and only if $\HH^{1}(A_{\lambda}) =0$.
\end{dem}

\begin{cor}
Let $A$ be an ada algebra. Then $A$ is simply connected if and only if the Hochschild cohomology ring is equal to $k$.
\end{cor}

\begin{dem}
This follows from \ref{theo:main} and \ref{cor:coho}(b).
\end{dem}


\begin{thebibliography}{10}

\bibitem{A} I.~Assem, {\it Left Sections and the left part of an artin algebra}, Colloquium Math.,{\bf  116} (2) (2009) 273-300.

\bibitem{ABCN} I.~Assem, J.C.~Bustamante, D.~Castonguay, C.~Novoa, {\it A note on the fundamental group of a one point extension}, Proyecciones {\bf 24} (1)(2005), 79-87.

\bibitem{AC1} I.~ Assem, F.~ U.~ Coelho, {\it Two-sided gluings of
tilted algebras}, J. Algebra {\bf 269} (2) (2003), 456-479.

\bibitem{AC2} I.~ Assem, F.~ U.~ Coelho, {\it Endomorphism algebras of projective modules modules over laura algebras}, J. Algebra and
Appl. {\bf 3} (1) (2004), 49-60.

\bibitem{ACT} I.~ Assem, F.~ U.~ Coelho, S.~ Trepode, {\it The left and the right parts of a module category}  J. Algebra  {\bf 281} (2) (2004), 518-534.

\bibitem{ACLST} I.~ Assem, F.U.~ Coelho, M.~ Lanzilotta, D.~Smith, S.~Trepode, {\it Algebras determined by their left and right parts}, Contemp. Math. {\bf 376}, Amer. Math. Soc., Providence, RI (2005) 13-47.

\bibitem{AL} I.~Assem, M. ~Lanzilotta, {\it The simple connectedness of a tame weakly shod algebra} Comm. Algebra {\bf 32} (9)(2004), 3685-3701.
\bibitem{ALR} I.~Assem, M. ~Lanzilotta, M. J. ~Redondo, {\it Laura Skew group algebras}, Comm. Algebra , {\bf 35} (7) (2007), 2241-2257.
\bibitem{AP} I.~Assem, J. A. ~de la Pe\~na, {\it The fundamental groups of a triangular algebra}, Comm. Algebra, {\bf 24} (1) (1996), 187-208.
\bibitem{ASS} I. ~Assem, D.~Simson, A.~ Skowro\'nski {\it Elements of the representation theory of associative algebras}, London Math. Soc. Student Texts {\bf 65}(2006) Cambridge Univ. Press, Cambridge.
\bibitem{AZ} I.~Assem, D.~Zacharia, {\it On split-by-nilpotent extensions}, Colloquium Math. {\bf 98}(2) (2003), 259-275.

\bibitem{ARS} M.~ Auslander, I.~ Reiten, S.~ Smal\o, {\it Representation theory of artin algebras}, Cambridge Studies in Advanced Mathematics {\bf 36} Cambridge University Press (1995) Cambridge.
\bibitem{AS} M.~ Auslander, S.~ Smal\o, {\it Almost split sequences in subcategories}. J. algebra {\bf 69} (1981) 426-454.

\bibitem{BG} K.~Bongartz, P.~Gabriel, {\it Covering spaces in representation-theory}, Invent. Math. {\bf 65} (3) (1981/82), 331-378.

\bibitem{C} F.U.~Coelho, {\it Directing components for quasitilted algebras}, Colloquium Math. {\bf 82} (1999) 271-275.

\bibitem{CL1} F.U.~Coelho, M.~A.~ Lanzilotta, {\it Algebras with small homological dimensions}, Manuscripta Math., {\bf 100}(1) (1999),1-11.

\bibitem{CL2}
F.~ U.~ Coelho, M.~ Lanzilotta, {\it Weakly shod algebras}, J. Algebra {\bf 265}(1) (2003), 379-403.


\bibitem{H1} D.~Happel, {\it Hochschild cohomology of finite dimensional algebras}. Sem. Marie-Paule Malliavin, Lect. Notes in Math. {\bf 1404}, Springer, Berlin (1989) 108-126.
\bibitem{H2} D.~Happel, {\it Triangulated categories in the representation theory of finite dimensional algebras}, London Math. Soc. Lecture Note Series {\bf 119}, Cambridge Univ. Press (1988).
\bibitem{H3} D.~Happel, {\it Hochschild cohomology of piecewise hereditary algebras}. Colloquium Math. {\bf 78} (1998) 261-266.
\bibitem{HRS} D.~Happel, I.~Reiten, S.~Smal\o, {\it Tilting in abelian categories and quasitilted algebras}, Proc. London Math. Soc.{\bf 46}(3) (1996).
\bibitem{L} P.~Le Meur, {\it Topological invariants of piecewise hereditary algebras}, Trans. Amer. Math. Soc. 363 (4) (2011), 2143-2170.

\bibitem{LS} H.~Lenzing, A.~Skowro\'nski, {\it Quasi-tilted algebras of canonical type}, Colloquium Math., {\bf 71} (2) (1996), 161-181.

\bibitem{RR} I.~Reiten, Ch.~Riedtmann, {\it Skew group algebras in the representation theory of Artin algebras}, J. Algebra, {\bf 92}, no.1, 224-282.
\bibitem{S1} A.~Skowro\'nski, {\it Simply connected algebras and Hochschild Cohomologies}, Proc. ICRA VI, Can. Math. Soc. Conf. Proc. {\bf 14} (1993) 431-447.

\bibitem{S2} A.~Skowro\'nski, {\it On artin algebras with almost all indecomposable modules of projective or injective dimension at most one}, Cent. Eur. J. Math. {\bf 1} (2003) 108-122.
\end{thebibliography}

\providecommand{\bysame}{\leavevmode\hbox
to3em{\hrulefill}\thinspace}

\end{document}